\documentclass[letterpaper, 10 pt, conference]{ieeeconf}  
 
\IEEEoverridecommandlockouts                              
\usepackage{graphics} 
\usepackage{graphicx}
\usepackage{amsmath}  
\usepackage{amssymb}  
\usepackage{amsthm} 
\usepackage{algpseudocode}
\usepackage{algorithm}
\usepackage{algorithmicx} 
\usepackage{setspace}
\usepackage{color}
\usepackage{verbatim} 
\usepackage[font= footnotesize]{caption}

\theoremstyle{plain} 
\newtheorem{assumption}{Assumption}
\theoremstyle{plain} 
\newtheorem{lem}{Lemma}
\theoremstyle{theorem} 
\newtheorem{thm}{Theorem}
\theoremstyle{remark}
\newtheorem{rem}{Remark}

\title{\LARGE \bf Asynchronous Splitting Design for Model Predictive Control{$^*$}}

\author{L. Ferranti$^{1}$, Y. Pu$^{2}$, C. N. Jones$^2$, and T. Keviczky$^1$
\thanks{*This research is supported by the European Union's Seventh Framework
Programme FP7/2007-2013 under grant agreement n. AAT-2012-RTD-2314544
(RECONFIGURE), by the TU Delft Space Institute, by the People Programme (Marie
Curie Actions) of the European Union's Seventh Framework Programme (FP7/2007-2013) under REA grant agreement n. 607957 (TEMPO) and by the Swiss National Science Foundation under Grand P2ELP$2\_165155$.} \thanks{{$^{1}$}L. Ferranti and T.
Keviczky are with the Delft Center for Systems and Control, Delft University of Technology, Delft, 2628 CD, The Netherlands,
        {\tt\small $\{$l.ferranti,t.keviczky$\}$@tudelft.nl}}%
\thanks{{$^{2}$}Y. Pu  and C. N. Jones are with the Laboratoire d'Automatique, \'{E}cole Polytecnique F\'{e}d\'{e}rale de Lausanne (EPFL), Lausanne, CH-1015, Switzerland,
        {\tt\small $\{$y.pu,colin.jones$\}$@epfl.ch}}}%

\newif\ifpaper 
\papertrue 
\begin{document}

\maketitle
\thispagestyle{empty}
\pagestyle{empty}

\begin{abstract}
This paper focuses on the design of an asynchronous dual solver suitable for embedded model predictive control (MPC) applications. 
The proposed solver relies on a state-of-the-art variance reduction (VR) scheme, previously used in the context of stochastic proximal gradient methods,
 and on the alternating minimization algorithm (AMA). The resultant algorithm, a stochastic AMA with VR, shows geometric convergence (in the expectation) to a suboptimal solution 
 of the MPC problem and, compared to other state-of-the-art dual asynchronous algorithms, allows to tune the probability of the asynchronous updates to improve the quality of the
  estimates. We apply the proposed algorithm to a specific class of splitting methods, i.e., the decomposition along the length of the prediction horizon, and provide preliminary 
  numerical results on a practical application, the longitudinal control of an Airbus passenger aircraft.
\end{abstract}
\section{INTRODUCTION}
\label{sec:sec_introduction}
Model Predictive Control (MPC) applications to systems with fast dynamics are still relatively limited~\cite{Keviczky2006,DiCairano2012}. Applications in 
fields such as automotive and aerospace have to deal very often with embedded legacy systems. These systems usually run on certified (for safety purposes) hardware architectures 
with limited availability, for example, of parallel computation units and support a small set of (certified) mathematical functions. In particular, the availability of optimization 
toolboxes suitable for MPC purposes on these platforms are limited (or nonexistent).  

Growing attention has been recently dedicated to the design of \emph{simple} first-order solvers for MPC~(\cite{Richter2009,Patrinos2014,StathopoulosECC13,Ferranti2015,Pu2014}). 
These solvers are relatively easy to certify (in terms of level of suboptimality of the solution), use only simple algebraic operations, and require little memory. 
In~\cite{StathopoulosECC13}-\cite{Pu2014}, operator-splitting methods, such as the alternating minimization method of multipliers (ADMM)~\cite{Boyd2011} and the fast alternating minimization 
algorithm (FAMA)~\cite{Goldstein2014}, have been used to exploit the MPC problem structure and speed-up the computation of the solution. These algorithms most of the time require frequent 
exchanges of information at given synchronization points. To reduce the bottlenecks at the synchronization points, a solver that can offer more~\emph{flexibility} in the way the solutions 
are computed (for example, by allowing asynchronous updates) would be attractive. In this work, we are interested in extending the use of splitting methods, such as AMA, to the asynchronous 
framework. 
    
{\it Contribution.}~The contribution of the paper is threefold. First, we propose a novel algorithm, a stochastic alternating minimization algorithm with variance reduction (SVR-AMA), 
suitable for MPC applications with state and input constraints. The proposed algorithm operates in the dual space and combines the advantages of the variance reduction scheme proposed 
in~\cite{Xiao2014} for the proximal stochastic gradient method with the alternating minimization algorithm~\cite{Tseng91}. The result is that the solution of the MPC problem can be computed 
in an asynchronous fashion (i.e., at each iteration, the algorithm updates a randomly selected subset of the dual variables instead of the whole set of dual variables) and the resultant 
algorithm has geometric convergence (in the expectation) to the optimal solution. Furthermore, the proposed algorithm allows the use of a generic probability distribution for the 
asynchronous updates. In addition, the probability distribution can be updated online to improve the quality of the estimates, as our numerical results show. Finally, the algorithm 
relies on simple algebraic operations, an appealing quality for embedded MPC applications.
  
Second, we show how we can use SVR-AMA for a specific splitting technique, i.e., the decomposition along the length of the prediction horizon 
(or \emph{time splitting}~\cite{StathopoulosECC13}).

Third, we present simulation results on a practical aerospace application, i.e., the longitudinal control of an Airbus passenger aircraft~\cite{Goupil2015}. The results
 show that the proposed method is more robust when solving ill-conditioned problems, outperforming synchronous methods in terms of computation time (measured in terms of number of iterations)
  and suboptimality level of the solution.
 
{\it Related Work.} SVR-AMA derives from the application to the dual problem of the proximal stochastic gradient method with variance reduction (Prox-SVRG) proposed in~\cite{Xiao2014}. 

\ifpaper The investigation of asynchronous dual algorithms for MPC is gaining more attention recently. In~\cite{Notarnicola2015}, for example, an asynchronous dual algorithm is proposed.
 Compared to~\cite{Notarnicola2015}, SVR-AMA allows the use of a generic (i.e., not necessarily uniform) probability distribution and, consequently, more flexibility in the tuning phase 
 of the algorithm.
\else \fi 

The idea of the time splitting has been previously proposed in~\cite{StathopoulosECC13}. Their work relies on a synchronous ADMM algorithm. In this context, we reformulate the approach
 for AMA to exploit SVR-AMA. 

\ifpaper
{\it Outline.} The paper is structured as follows. Section~\ref{sec:problem_formulation} introduces the MPC problem formulation. Section~\ref{sec:sec_preliminaries} summarizes AMA and 
Prox-SVRG. Section~\ref{sec:sec_async_algorithm} describes SVR-AMA and show convergence results.
Section~\ref{sec:async_mpc} shows how to reformulate the proposed MPC problem for SVR-AMA using the time splitting. Section~\ref{sec:sec_simulation_results} presents
 preliminary numerical results using an aerospace example. Finally, Section~\ref{sec:sec_conclusions} concludes the paper.
\fi
 
{\it Notation.} For $u\in \mathbb R^n$, $\left \| u \right
\|=\sqrt{\langle u,u \rangle}$ is the Euclidean norm. Let $\mathbb C$ be a convex set. Then, $\mathbf{Pr}_{\mathbb C}(u)$ is the projection of $u$ onto $\mathbb C$. 
Let $f:\mathcal D \rightarrow \mathcal C$ be a function. Then, $f^{\star}(y)=\operatorname{sup}_x(y^{\operatorname{T}}x-f(x))$ and $\nabla f(x)$ are the conjugate 
function and the gradient of $f(x)$, respectively. Furthermore, $\mathcal I_{\mathbb C}(\sigma)$ is the indicator function on the convex set $\mathbb C$, which
 is zero if $\sigma \in \mathbb C$ and infinity otherwise.
Let $A\in\mathbb R^{n\times m}$. Then, $\operatorname{eig}_{\max}(A)$ and $\operatorname{eig}_{\min}(A)$ are the largest and the smallest (modulus) eigenvalues of $A^{\textrm{T}}A$. $P  \in\mathbb S_{+}^{n\times n}$ denotes that $P\in \mathbb R^{n\times n}$ is positive definite. In addition, let $x\in \mathbb R^n$ be a random variable, $\mathbb E [x]$ is its expected value. Finally, details on the notions of strong convexity and Lipschitz continuity used in the paper can be found in~\cite{Boyd2004}.
 
\section{PROBLEM FORMULATION}
\label{sec:problem_formulation}
Consider the discrete linear time-invariant (LTI) system described by the following equation:
\begin{equation}
\label{eq:state_space}
x(k+1) = A x(k) + B u(k),\quad k=0,1,2,\ldots
\end{equation}
The state $x(k)\in \mathbb R^n$ and the control input $u(k)\in\mathbb R^m$ are subject to the following polyhedral constraints:
\begin{equation}
\label{eq:constraints}
Cx(k) + Du(k)\leq d,
\end{equation}
where $C\in\mathbb R^{p\times n}$ and $D \in \mathbb R^{p\times m}$. Note that the definition of the constraints~\eqref{eq:constraints} 
can include constraints on $x(k)$ only or on $u(k)$ only. 
We aim to regulate the state $x(k)$ to the origin using the control input $u(k)$ while respecting the constraints~\eqref{eq:constraints}. 
This goal can be translated into the following model predictive control (MPC) problem:
\begin{subequations}
\label{eq:initial_MPC_problem}
\begin{align}
&\min\limits_{x,u}~\frac{1}{2}\sum\limits_{t = 0}^{N} \left(x_t^{\textrm{T}}Qx_t+u_t^\textrm{T}Ru_t\right)\\
\label{eq:initial_MPC_problem_dyn_coupling}
&\textrm{s.t.:~} x_{t+1} = Ax_t+Bu_t \quad t=0,\ldots,N-1\\
&~\quad~ Cx_t+Du_t\leq d \quad\quad~ t=0,\ldots,N\\
&~\quad~ x_0 = x_{\textrm{init}},
\end{align}
\end{subequations}
where $x_t$ and $u_t$ represent the $t$-step-ahead state and control predictions, respectively, $N$ indicates the length of the prediction horizon, 
$Q\in \mathbb S_{+}^{n\times n}$, $R\in \mathbb S_{+}^{m\times m}$, and $x_\textrm{init}$ is the initial (measured) state vector. The MPC
 law implemented in closed loop is given by the first element of the optimal control sequence obtained by solving Problem~\eqref{eq:initial_MPC_problem}, i.e., $u_{\textrm{MPC}} = u^*_0$.

Our goal is to solve Problem~\eqref{eq:initial_MPC_problem} in an embedded environment. In particular, we assume that explicit MPC~\cite{Bemporad2002} cannot 
be used due to the problem size and that the computational resources are limited\ifpaper, i.e., parallel architectures are not available, memory resources are limited, and only simple 
algebraic operations are supported.\fi~With this framework in mind, in the following, we focus on the design of a simple solver for Problem~\eqref{eq:initial_MPC_problem} that 
relies on operator-splitting methods\ifpaper~(which, for example, usually rely on parallel hardware architectures)\fi and asynchronicity (which allows one to perform updates of
 a randomly selected subset of variables to reduce the computational effort). The next section introduces the techniques we rely on to solve Problem~\eqref{eq:initial_MPC_problem}, 
 i.e., AMA proposed by~\cite{Tseng91} and the proximal stochastic gradient descent method with variance reduction (Prox-SVRG) proposed by~\cite{Xiao2014}.

\section{PRELIMINARIES}
\label{sec:sec_preliminaries}
AMA and Prox-SVRG are the main techniques we rely on for our design. In the following, Section~\ref{subsec:ama} briefly describes
 AMA developed by~\cite{Tseng91}, while Section~\ref{subsec:prox_svrg} describes Prox-SVRG~\cite{Xiao2014}. 
\subsection{Alternating Minimization Algorithm}
\label{subsec:ama}
\captionsetup[algorithm]{font=normal}
\begin{algorithm}[t]
\begin{algorithmic}  
 \State{Given $\mu^{\textrm{0}}$, $T$, and $\tau < \sigma_f/\operatorname{eig}_{\max}(H_y)$.}
\While{$k = 1, \ldots, T$}
\State{{1a.} $y^k=\operatorname{argmin}_y~ f(y)+\langle\mu^{k-1},-H_yy\rangle$.}
\vspace{0.02in}
\State{{1b.} $z^k=\operatorname{argmin}_z~ g(z)+\langle\mu^{k-1},-H_zz\rangle+$}
\State{$~~~~~~~~~~~~~~~~~~~~~+\frac{\tau}{2}\|d-H_yy^k-H_zz\|^2$.}
\vspace{0.02in}
\State{{2.} $\mu^k = \mu^{k-1} +\tau (d-H_yy^{k}-H_zz^K)$}
\EndWhile
 \caption{AMA~\cite{Tseng91}.}
 \label{alg:ama}
 \end{algorithmic} 
\end{algorithm}
Consider the following problem:
\begin{subequations}
\label{eq:ama_problem}
\begin{align}
\operatorname{minimize}~ &f(y) + g(z)\\
\label{eq:ama_eq_constraints}
\operatorname{subject~to}~& H_yy +H_z z = d,
\end{align}  
\end{subequations} 
where $f(y) := \sum_{t=0}^Nf_t(y)$ under the following assumptions:
\begin{assumption}
\label{ass:strong_convexity_f}
$f_t$ is a strongly convex function and $\sigma_{f_t}$ denotes its convexity parameter ($t=0,\ldots,N$).
\end{assumption}
\begin{assumption}
\label{ass:lipschitz_continuity_f}
$f_t$ has a Lipschitz continuous gradient with modulus $L_{f_t}$ ($t=0,\ldots,N$) and $\operatorname{eig}_{\min}(H_y)>0$.
\end{assumption}
\begin{assumption}
\label{ass:convexity_g}
$g$ is a convex function not necessarily smooth. 
\end{assumption}
\ifpaper
Furthermore, recall the following properties of the conjugate function $f^\star$:
\begin{lem}[Thm. 4.2.1~\cite{Hul93}]
\label{lem:property_strong_convexity}
If $f$ is strongly convex with convexity parameter $\sigma_f$, then $f^{\star}$ has a Lipschitz continuous gradient with modulus $L(\nabla f^{\star}) = \sigma_f^{-1}$.
\end{lem}
\begin{lem}[Thm. 4.2.2~\cite{Hul93}]
\label{lem:property_lipschitz_continuity}
If $f$ is convex and has a Lipschitz continuous gradient with modulus
$L_f$, then $f^\star$ is strongly convex with convexity parameter $L_f^{-1}$.
\end{lem}
\else 
\fi   
A state-of-the-art algorithm to solve Problem~\eqref{eq:ama_problem} is AMA~\cite{Tseng91}. AMA operates as a proximal gradient algorithm (such as, ISTA~\cite{Beck}) on the dual of Problem~\eqref{eq:ama_problem}. Specifically, given the dual of Problem~\eqref{eq:ama_problem} (under the assumptions above), described as follows:
\begin{equation}
\label{eq:ama_dual_problem}
\underset{\mu\in\mathbb R^{n_{\mu}}}{\operatorname{maximize}}~ D(\mu)\left\{:= -F(\mu)-G(\mu)\right\}, 
\end{equation}
where $F(\mu) := \sum_{t=0}^Nf_t^{\star}(H_{y_t}^{\textrm{T}}\mu)$ and $G(\mu):= g^{\star}(H_z^{\textrm{T}}\mu)-d^{\textrm{T}}\mu$, the following holds:
\begin{lem}
\label{lem:strong_convexity_F}
If Assumptions~\ref{ass:strong_convexity_f}-\ref{ass:convexity_g} are satisfied, $F(\mu)$ is strongly convex with Lipschitz continuous gradient characterized by Lipschitz constant $L(\nabla F):= \operatorname{eig}_{\max}(H_y)\sigma_f^{-1}$. Furthermore $G(\mu)$ is convex.
\end{lem}
\ifpaper
\begin{proof}
We can use Lemmas~\ref{lem:property_strong_convexity} and \ref{lem:property_lipschitz_continuity} to derive the properties of $F(\mu)$. Convexity of $G(\mu)$ follows from the properties of the conjugate of a convex function and from the fact that $d^{\textrm{T}}\mu$ is a linear function.  
\end{proof}
\fi
AMA updates the dual variables $\mu\in\mathbb R^{p_{\mu}}$ as described in Algorithm~\ref{alg:ama}. In general, AMA uses only simple algebraic operations (if $y$ and $z$ are unconstrained steps 1a and 1b can be performed efficiently) and does not require advanced hardware architectures. Nevertheless, the algorithm requires frequent exchange of information at given synchronization points (e.g., step 1b requires $y$ computed at step 1a that can lead to bottlenecks in the computation of the problem solution). Hence, it would be better to have some flexibility on the update strategy. Motivated by this observation, the following section introduces Prox-SVRG used to derive our proposed asynchronous AMA, as described in Section~\ifpaper\ref{sec:sec_async_algorithm}\fi.   
\subsection{Stochastic Gradient Method with Variance Reduction} 
\label{subsec:prox_svrg}
 Consider the following primal problem:
 \begin{equation} 
 \label{eq:svrg_problem}
 \underset{y\in\mathbb R^{n_{y}}}{\operatorname{minimize}}~P(y)\left\{:=F(y)+G(y)\right\},
 \end{equation}
 where $F(y):=\sum_{t=0}^{N}F_t(y)$ and $G(y)$ satisfy the following assumptions: 
\begin{assumption}
\label{ass:strong_convexity_F}
$F(y)$ is a strongly convex function with convexity parameter $\sigma_F$ and Lipschitz continuous gradient characterized by a Lipschitz constant $L\leq \sum_{t=0}^N L_t$, where $L_t$ are the Lipschitz constants of each $F_t (y)$.
\end{assumption}
\begin{assumption}
\label{ass:convexity_G}
$G(y)$ is a convex function.
\end{assumption} 
Furthermore, define the~\emph{proximal operator} as follows:
\begin{equation}
\label{eq:svrg_proximal_operator}
\textrm{prox}_{\tau G} (y) := \underset{y\in\mathbb{R}^{n_{y}}}{\operatorname{argmin}}\left\{\frac{1}{2}\left\|y-x\right\|^2+\tau G(y)\right\}.
\end{equation}  
 
The main idea behind Prox-SVRG~\cite{Xiao2014} is to eliminate the dependency of the number of iterations (typical of stochastic gradient methods) in the definition of the step size and reduce the burden in the computation of $\nabla F(y)$ (typical of classical gradient methods). As pointed out in~\cite{Xiao2014}, proximal stochastic gradient methods (such as,~\cite{Duchi2009,Langford2009}) suffer of sublinear convergence (to a suboptimal solution of Problem~\eqref{eq:svrg_problem}) given that the step size decreases at each iteration of the algorithm, but behave well when $N$ is large. On the other hand, classical proximal gradient methods require at each iteration of the algorithm to compute the full gradient of $F(y)$, which can be an involved operation if $N$ is large, but the step size is fixed and independent of the number of iterations (leading to better theoretical convergence properties). Hence, Prox-SVRG aims to exploit the benefits of the two techniques as explained below and described in Algorithm~\ref{alg:prox_svrg}.   

 \captionsetup[algorithm]{font=normal}
\begin{algorithm}[t]
\begin{algorithmic}  
 \State{Given $\tilde y^{\textrm{0}}$, $N$, $\mathcal I_N:=\{0,\ldots,N\}$ $\tau$, and $T$.}
\While{$s\le\bar s$}
\State{{0a.} Set $\tilde y = \tilde y^{s-1}$.}
\State{{0b.} Set $\tilde \beta = \nabla F(\tilde y)$.}
\State{{0c.} Set $y^0 = \tilde y$.}
\State{{0d.} Set $\Pi:=\{\pi_0,\ldots,\pi_N\}$.}
\For{$k = 1, \ldots, T$}
\State{{1.} Pick $i\in\mathcal I_N$ randomly according to $\Pi$.}
\State{{2.} $\beta^k = \tilde \beta+\frac{\nabla F_{i}(y^{k-1})-\nabla F_{i}\left(\tilde y\right)}{\pi_{i}}$}
\State{{3.} $y^k = \textrm{prox}_{\tau G} (y^{k-1}-\tau \beta^k)$}
\EndFor
\State{{4.} $\tilde y^s = \sum_{k=1}^T y^k$.}
 \EndWhile
 \caption{Prox-SVRG~\cite{Xiao2014}.}
 \label{alg:prox_svrg}
 \end{algorithmic} 
\end{algorithm}

Prox-SVRG uses a multistage strategy to gradually reduce the variance in the estimation of the full gradient $\nabla F(y)$ (without computing the actual full gradient at each iteration). In particular, the full gradient of $F(y)$ is updated only every $T$ iterations to reduce the computational effort compared to the classical gradient methods, and the proximal step (step~3) uses a modified direction $\beta^k$ (step 2) that leads to a smaller variance $\mathbb E\|\beta^k-\nabla F(y^{k-1})\|^2$ compared to the one obtained using classical stochastic gradient methods $\mathbb E\|\nabla F_{i}(y^{k-1})-\nabla F(y^{k-1})\|^2$ ($i\in\mathcal I_{N}$), where $\nabla F_{i}(y^{k-1})$ is used as update direction (refer to~\cite{Xiao2014} for more details). Furthermore, the random sampling (step 1) is performed on a probability distribution $\Pi:=\{\pi_0,\ldots,\pi_N\}$ that does not necessarily have to be uniform, i.e., the algorithm allows more flexibility in the tuning phase by supporting other distributions as well, such as Poisson distributions, normal distributions, etc. Algorithm~\ref{alg:prox_svrg} achieves geometric convergence in the expectation, as stated in the following theorem:
\begin{thm}[{Thm. 3.1} in~\cite{Xiao2014}]
\label{th:prox_svrg} 
Suppose Assumptions~\ref{ass:strong_convexity_F} and~\ref{ass:convexity_G} hold. Let $y^{*}~=~{\operatorname{argmin}}_y~P(y)$ and $L_{\Pi} := \max_t L_t/\pi_t$. Assume that $0<\tau<1/(4L_{\Pi})$ and $T$ is sufficiently large so that:
\begin{equation}
\rho :=\frac{1}{\tau\sigma_FT(1-4\tau L_{\Pi})}+\frac{4 \tau L_{\Pi}(T+1)}{T(1-4\tau L_{\Pi})}<1.
\end{equation}
Then, for $\bar s >1$, Algorithm~\ref{alg:prox_svrg} has geometric convergence in expectation:
\begin{equation}
\mathbb E P(\tilde y^{\bar s})-P(y^*)\leq \rho^{\bar s} \left[P(\tilde y^0)-P(y^*)\right].
\end{equation}
\end{thm}
\begin{rem}
\label{rem:probability}
The dependency on the probability $\pi_t$ in the choice of the step size $\tau$ can be problematic when using probability distributions with $\pi_t\to 0$ or when $N\to\infty$ (e.g., the constrained infinite horizon LQR). Nevertheless, this dependency can be removed in the special case in which $F(y) = \sum_{t=1}^{N}F_t(y_i)$, i.e., when the cost is separable in $y_t$. 
From the MPC perspective, this is very often the case when the dual formulation is used, as it is shown in Section~\ref{sec:async_mpc} for the decomposition along the length of the prediction horizon. Hence, this observation, when the algorithm is used in the dual framework, can be very beneficial to improve the choice of the step size and the quality of the MPC solution. 
\end{rem}
\ifpaper \section{STOCHASTIC AMA WITH VARIANCE REDUCTION}
\label{sec:sec_async_algorithm}
\fi
\ifpaper
Our goal is to solve Problem~\eqref{eq:ama_problem} in an asynchronous fashion, i.e. by allowing updates of a randomly selected subset
 of the dual variables at each iteration of the solver. Hence, given that Algorithm~\ref{alg:prox_svrg} cannot be directly applied to 
 Problem~\eqref{eq:ama_problem}, we proceed as explained in Section~\ref{subsec:ama}, i.e., we apply Algorithm~\ref{alg:prox_svrg} to the dual of 
 Problem~\eqref{eq:ama_problem}. The resultant algorithm (SVR-AMA) is described by Algorithm~\ref{alg:svrama}.
\fi
\begin{thm}
\label{th:svg_ama}  
Suppose Assumptions~\ref{ass:strong_convexity_f}-\ref{ass:convexity_g} hold. Let $\mu^{*}~=~{\operatorname{argmax}}_{\mu}~D(\mu)$, where $D(\mu)$ 
is the dual cost defined in~\eqref{eq:ama_dual_problem}. In addition, let $L_{\Pi}^{\star} := \max_t (\pi_t\sigma_{f})^{-1}\operatorname{eig}_{\min}(H_y)$, $\pi_t\in\Pi$.
 Assume that $0<\tau<1/(4L_{\Pi}^{\star})$ and $T\geq 1$ such that:
\begin{equation}
\rho :=\frac{L_f}{\tau T(1-4\tau L_{\Pi}^\star)}+\frac{4 \tau L_{\Pi}^\star(T+1)}{T(1-4\tau L_{\Pi}^\star)}<1.
\end{equation}
Then, for $\bar s>0$, Algorithm~\ref{alg:svrama} has geometric convergence in expectation:
\begin{equation}
 D(\mu^*) - \mathbb E D(\tilde \mu^{\bar s})\leq \rho^{\bar s} \left[D(\mu^{*})-D(\tilde \mu^0)\right].
\end{equation}
\end{thm} 
\ifpaper
\begin{proof}
If the assumptions of the theorem are satisfied, exploiting a similar argument to the one in Theorem 1 in~\cite{Goldstein2014} for AMA, 
we can show that Algorithm~\ref{alg:svrama} is equivalent to applying Algorithm~\ref{alg:prox_svrg} to Problem~\eqref{eq:ama_dual_problem}. Consequently, by using 
Lemmas~\ref{lem:property_strong_convexity} and~\ref{lem:property_lipschitz_continuity}, the results follows from Theorem~\ref{th:prox_svrg}.  
\end{proof}
\fi
\ifpaper 
\begin{rem}
\label{rem:svrfama}
Note that an acceleration step such as the one used in FAMA~\cite{Goldstein2014} can be added to improve the convergence of the algorithm. 
The proof of convergence for the acceleration is part of our future work, but note that the simulation results in Section~\ref{sec:sec_simulation_results} rely on 
the acceleration step of FAMA, showing that, in practice, it can be applied to the proposed algorithm. In addition, note that an acceleration strategy for Prox-SVRG~\cite{Xiao2014} 
has been recently proposed in~\cite{Nitanda2014} and can be extended to the dual framework to accelerate Algorithm~\ref{alg:svrama}. 
\end{rem}
\else
\fi

\ifpaper\section{ASYNCHRONOUS MPC}
\label{sec:async_mpc}
\else
\fi
\ifpaper Our aim is to solve the MPC Problem~\eqref{eq:initial_MPC_problem} presented in Section~\ref{sec:problem_formulation} using Algorithm~\ref{alg:svrama}.
 Hence, we must show that the MPC Problem~\eqref{eq:initial_MPC_problem} is a particular case of Problem~\eqref{eq:ama_problem}.
\fi 

First, we decompose Problem~\eqref{eq:initial_MPC_problem} along the length of the prediction horizon $N$ into $N+1$ smaller subproblems, according to the \emph{time-splitting} 
strategy proposed in~\cite{StathopoulosECC13}. This results is achieved thanks to the introduction of $N$ consensus variables $z_t\in \mathbb R^n$ ($t=1,\ldots,N$) used 
to break up the dynamic coupling~\eqref{eq:initial_MPC_problem_dyn_coupling}. This decomposition allows to reformulate Problem~\eqref{eq:initial_MPC_problem} as follows:
\begin{subequations}
\label{eq:time_splitting_MPC_problem}
\begin{align}
&\min\limits_{x,u}~\frac{1}{2}\sum\limits_{t = 0}^{N} \left(x_t^{(t)^{\textrm{T}}}Qx_t^{(t)}+u_t^{(t)^{\textrm{T}}}R u_t^{(t)}\right)\\
\label{eq:time_splitting_MPC_problem_consensus_1}
&\textrm{s.t.:~}z_{t+1} = Ax_t^{(t)}+Bu_t^{(t)} \quad t=0,\ldots,N-1\\
\label{eq:time_splitting_MPC_problem_consensus_2}
&~\quad~ z_{t+1} = x_{t+1}^{(t+1)} \quad\quad\quad\quad~ t=0,\ldots,N-1\\
\label{eq:time_splitting_MPC_problem_ineq_constraints}
&~\quad~ Cx_t^{(t)}+Du_t^{(t)}\leq d \quad\quad~ t=0,\ldots,N\\
&~\quad~ x_0^{(0)} = x_{\textrm{init}},
\end{align}
\end{subequations}
where the original dynamic coupling~\eqref{eq:initial_MPC_problem_dyn_coupling} has been replaced by the consensus constraints~\eqref{eq:time_splitting_MPC_problem_consensus_1}
 and~\eqref{eq:time_splitting_MPC_problem_consensus_2}. Note that we introduced the superscript $t$ to underline that the $x_t$ and $u_t$ are local variables of the subproblems 
 obtained after the time splitting~\cite{StathopoulosECC13}.
Finally, if we introduce $N+1$ additional slack variables $\sigma_t\in\mathbb R^{p}$ to remove the inequality constraints~\eqref{eq:time_splitting_MPC_problem_ineq_constraints} 
and define $\mathbb C:=\left\{\sigma_t\in \mathbb R^{p}\,|\,\sigma_t\geq 0\right\}$,  Problem~\eqref{eq:time_splitting_MPC_problem} can be written as the sum of the following 
subproblems:
\captionsetup[algorithm]{font=normal}
\begin{algorithm}[t]
\begin{algorithmic}  
 \State{Given $\tilde \mu^{\textrm{0}}$, $N$, $\mathcal I_N:=\{0,\ldots,N\}$, $\tau$, and $T$.}
\While{$s\le\bar s$}
\State{{0a.} Set $\tilde \mu = \tilde \mu^{s-1}$, $\tilde y = \tilde y^{s-1}$.}
\State{{0b.} Set $\tilde \beta = \nabla F(\tilde \mu)$.}
\State{{0c.} Set $\mu^0 = \tilde \mu$.}
\State{{0d.} Set $\Pi:=\{\pi_0,\ldots,\pi_N\}$.}
\For{$k = 1, \ldots, T$}
\State{{1.} Pick $i\in\mathcal I_N$ randomly according to $\Pi$.}
\State{{2a.} $y^k=\operatorname{argmin}_y~ f_{i}(y)+\langle\mu^{k-1},-H_{y_{i}}y\rangle$.}
\State{{2b.} $z^k=\operatorname{argmin}_z~ g(z)+\langle\mu^{k-1},-H_{z}z\rangle+$\\$~~~~~~~~~~~~~~~~~~~~~~~~~~~~~~~+\frac{\tau}{2}\|d-H_yy^k-H_zz\|^2$.}
\vspace{0.05in}
\State{{3.} $\beta^k = \tilde \beta+{(y^k-\tilde y)^{\textrm{T}}H_{y_i}^{\textrm{T}}}{\pi_{i}^{-1}}$}
\State{{4.} $\mu^k = \mu^{k-1}-\tau (\beta^k+H_{z}z^k-d)$}
\EndFor
\State{{5.} $\tilde \mu^s = \frac{1}{T}\sum_{k=1}^T \mu^k$, $\tilde y^s = \frac{1}{T}\sum_{k=1}^T y^k$.}
 \EndWhile
 \caption{SVR-AMA.}
 \label{alg:svrama}
 \end{algorithmic} 
\end{algorithm}
 \begin{subequations}
				  \label{eq:subproblem_t_primal}
					\begin{align}
					\label{eq:subproblem_t_primal_cost}
                     &\min\limits_{y_t} f_t\left(y_t\right)+\sum_{i=1}^{p} \mathcal I_{\mathbb                         C_{}}(\sigma_{t_i})\\
                     \label{eq:subproblem_t_primal_w}
                     &\textrm{s.t.:~}~~~w_t:~~z_{t_{~}} = H_1 y_{t},\\
                     \label{eq:subproblem_t_primal_v}
                     &\quad~~~v_{t+1}:~z_{t+1} = H_2 y_{t},\\
                     &\quad~~~~~~ \lambda_t:~~~\sigma_{t}=d-G y_{t},
					\end{align}
\end{subequations}
where we define $y_t^{\textrm{T}}:=[x_{t}^{(t)^{\textrm{T}}}~u_{t}^{(t)^{\textrm{T}}}]$, 
$f_t(y_t):=y_t^{\textrm{T}}\mathcal Qy_t$, $\mathcal Q :=\textrm{diag}\left\{Q,R\right\}$, $G := \left[C~D\right]$, $H_1:=\left[I_n~0_{n\times m}\right]$, $H_2:= \left[A~B\right]$. 
Furthermore, for each equality constraint in Problem~\eqref{eq:subproblem_t_primal}, the corresponding Lagrange multipliers have been highlighted. 

If we define $\mathbf y^{\operatorname{T}} = [y_0^{\operatorname{T}}\ldots y_N^{\operatorname{T}}]$, $f(\mathbf y) = \mathbf{y}^{\operatorname{T}}\mathbf Q \mathbf y=\sum_{t=0}^nf_t(y_t)$, 
$\mathbf z^{\operatorname{T}}~=~[z_1^{\operatorname{T}}\ldots z_N^{\operatorname{T}}~ \sigma_0^{\operatorname{T}}\ldots\sigma_N^{\operatorname{T}}]$, 
$\mathbf Q~=~\operatorname{diag}\{\mathcal Q\ldots\mathcal Q\}$, $g(\mathbf{z})=\sum_{t=0}^N\mathcal{I}_{\mathbb C}(\sigma_t)$, $h_{y_0}^{\operatorname{T}} = 
\left[H_2^{\operatorname{T}}~|-G^{\operatorname{T}}\right]$, $h_{y}^{\operatorname{T}} = \left[H_1^{\operatorname{T}}~ H_2^{\operatorname{T}}~| -G^{\operatorname{T}}\right]$,
 $h_{y_N}^{\operatorname{T}} = \left[H_1^{\operatorname{T}}~| -G^{\operatorname{T}}\right]$, $h_{d_0}^{\operatorname{T}}=$ $ h_{d_N}^{\operatorname{T}}  =\left [0_n~ |-d^{{\operatorname{T}}}\right]$, $h_{d}^{\operatorname{T}} = \left[0_n ~0_n~ |-d^{{\operatorname{T}}}\right]$,
\begin{align*}
H_{\mathbf y}& := \begin{bmatrix}h_{y_0}    & 0   & \ldots & \ldots & 0\\
                                  0      & h_y &   & &0\\
                                  \vdots &    &\ddots   &  & \vdots\\
                                  0      & 0   & \ldots  & h_{y} & 0\\
                                  0      & 0   & \ldots  & 0& h_{y_N}
\end{bmatrix},~\mathbf{d} :=\begin{bmatrix}
 h_{d_0} \\h_d\\\vdots\\h_d\\h_{d_N}
\end{bmatrix},\\
H_{\mathbf z}& := \left[\begin{array}{cccc|cccc}
I_n & 0 &\ldots& 0 & 0 & 0 &\ldots & 0\\
0   & 0 &\ldots& 0 & I_p & 0 &\ldots & 0\\
\hline
I_n & 0 & \ldots & 0 & 0 & 0 & \ldots & 0\\
0 & I_n & \ldots & 0 & 0 & 0 & \ldots & 0\\
0 & 0 & \ldots & 0 & 0 & I_p& \ldots& 0\\
\hline
\vdots & & \ddots &  & \vdots &  & \ddots & 0\\
\hline
0 & 0 & \ldots & I_n & 0 & 0 & \ldots & 0\\
0 & 0 & \ldots & 0 & 0 & 0 & \ldots & I_p\\
\end{array}\right],
\end{align*}
Problem~\eqref{eq:initial_MPC_problem} can be rewritten as follows:
\begin{subequations}
\label{eq:splitting_problem_condensed}
\begin{align}
{\operatorname{minimize}}~f(\mathbf{y}) + g(\mathbf{z})\\
H_{\mathbf y} \mathbf y + H_{\mathbf{z}} \mathbf z= \mathbf d, 
\end{align}
\end{subequations}
According to the definition of $\mathcal Q$, $f(\mathbf y)$ is strongly convex, has a convexity modulus 
$\sigma_f:=\textrm{eig}_{\min}(\mathcal Q)= \operatorname{eig}_{\min}(\operatorname{blockdiag}\{Q,R\}) = \sigma_{f_t}$, and a Lipschitz constant
 $L_f:= \operatorname{eig}_{\max}(\mathcal Q)$. In addition, $g(\mathbf z)$ is a convex function. 

For the proposed splitting, Assumptions~\ref{ass:lipschitz_continuity_f} and~\ref{ass:convexity_g} are satisfied. Concerning Assumption~\ref{ass:strong_convexity_f}, 
note that 
$L(\nabla F):= \operatorname{eig}_{\max}(H_{\mathbf y})\sigma_{f}^{-1}= \max_t(\operatorname{eig}_{\max}(H_{y_t})\sigma_{f_t}^{-1}) = 
\max_t(L_t(\nabla F_t)) = L_t(\nabla F_t) \leq \sum_{t=0}^N L_t(\nabla F_t)$, where the last equality follows from the fact that we deal with LTI systems
 ($L_0 = L_1 =\ldots= L_N$). Hence, on the dual, Assumption~\ref{ass:strong_convexity_F} still holds and, consequently, we can use SVR-AMA to solve 
 Problem~\eqref{eq:initial_MPC_problem}. 

\ifpaper
The associated SVR-AMA algorithm to solve Problem~\eqref{eq:splitting_problem_condensed} is detailed in Algorithm~\ref{alg:splittingAMA}.
In particular, defining 
$\boldsymbol \mu^{\textrm{T}} :=[v_1^{\textrm{T}}~\lambda_0^{\textrm{T}}~|~w_1^{\textrm{T}}~\ldots~w_{N-1}^{\textrm{T}}~v_{N}^{\textrm{T}}~\lambda_{N-1}^{\textrm{T}}~|
~w_N^{\textrm{T}}~\lambda_{N}^{\textrm{T}}]$, according to the partitioning of 
$H_{\mathbf y}$ and $H_{\mathbf{z}}$, $F(\boldsymbol \mu)= f^{\star}(H_{\mathbf y}^{\textrm{T}}\boldsymbol \mu)$. 
Furthermore, $\nabla F(\mathbf w)$ is the gradient of $F$ at $\mathbf w$, $\nabla F(\mathbf v)$ 
is the gradient of $F$ at $\mathbf v$, and $\nabla F(\boldsymbol \lambda)$ is the gradient of $F$ at $\boldsymbol \lambda$. 
Note that the calculation of the gradient step for this particular splitting is very simple and requires the evaluation of the product $H_{\mathbf y}\mathbf y$, 
that can be performed efficiently by exploiting the structure of the matrix $H_{\mathbf y}$. Finally, note that, given the structure of $F(\boldsymbol \mu)$ 
the probability $\pi_t$ does not affect the choice of the step size $\tau$, according to Remark~\ref{rem:probability}.
 
The following complexity upper bound on the primal sequence can be defined:
\begin{thm}
Consider Problem~\eqref{eq:splitting_problem_condensed}. Let $\{\mathbf y^k\}$ and $\{\boldsymbol \mu^k\}$ be the sequence 
of primal and dual variables, respectively, generated by Algorithm~\ref{alg:splittingAMA}.
 If Assumptions~\ref{ass:strong_convexity_f}-\ref{ass:convexity_g} are satisfied, given  $\tilde {\boldsymbol \mu}^0 \in \operatorname{dom}(G)$, 
 where $G:= g^{\star}(H_{\mathbf z}^{\textrm{T}}\boldsymbol \mu)-\boldsymbol d^{\textrm{T}}\boldsymbol\mu$, then, the following holds:
\begin{equation}
\label{eq:th3}
\mathbb E\|\tilde {\boldsymbol y}^s-\boldsymbol y^*\|^2\leq \frac{2}{\sigma_f}(D(\boldsymbol \mu^*) -\mathbb E D(\tilde{\boldsymbol \mu})).
\end{equation}
\end{thm} 
\captionsetup[algorithm]{font=normal}
\begin{algorithm}[t]
\begin{algorithmic}  
 \State{Given $\tilde \mu^{\textrm{0}}$, $N$, $\mathcal I_N:=\{0,\ldots,N\}$, $L^{\star} := (\sigma_{f})^{-1}\operatorname{eig}_{\max}(H_y)$, $0<\tau<1/(4L^{\star})$, and $T$.}
\While{$s\le\bar s$}
\State{{0a.} Set $\tilde {\mathbf w} = \tilde {\mathbf w}^{s-1}$, $\tilde {\mathbf v} = \tilde {\mathbf v}^{s-1}$, \hfill \\
~~~~~~~~ $\tilde {\boldsymbol\lambda} = \tilde {\boldsymbol\lambda}^{s-1}$, and $\tilde {\boldsymbol y} = \tilde {\boldsymbol y}^{s-1}$.}
\State{{0b.} Set $\tilde \beta_{\mathbf w} = \nabla F(\tilde {\mathbf w})$, $\tilde \beta_{\mathbf v} = \nabla F(\tilde {\mathbf v})$, and\hfill\\~~~~~~~~~~~~ $\tilde \beta_{\boldsymbol \lambda} = \nabla F(\tilde {\boldsymbol \lambda})$.}
\State{{0c.} Set $\mathbf {w}^0 = \tilde {\mathbf w}$, ${\mathbf v}^0 = \tilde {\mathbf v}$, and $\boldsymbol \lambda^0 = \tilde {\boldsymbol\lambda}$.}
\State{{0d.} Set $\Pi:=\{\pi_0,\ldots,\pi_N\}$ on $\mathcal I_N$}
\For{$k = 1, \ldots, T$}
\State{{1.} Pick $i\in\mathcal I_N$ randomly according to $\Pi$.}
\State{{2a.} $y_{i}^k=\operatorname{argmin}_y~ f_{i}(y_{i})+\langle w_{i},H_1y_{i}\rangle+$\hfill\\ ~~~~~~~~~~~~~~~~~~~~~$\langle v_{i+1},H_2y_{i}\rangle+\langle\lambda_{i},-Gy_{i}\rangle$.}
\State{{2.b} $\sigma_{i}^{k}=\mathbf{Pr}_{\mathbb C_{}}(G y_{i}^{k}-d-\frac{1}{\tau} \lambda_{i})$.}
\vspace{0.05in}
\State{{2c.} $z_{i}^{k} = \frac{1}{2}\left[{H_1 y_{i}^{k}+H_2 y_{i-1}}-\frac{1}{\tau}(w_{i}+ v_{i})\right]$.}
\vspace{0.05in}
\State{{3a.} $\beta_{w_{i}}^k = \tilde \beta_{w_{i}}+\frac{(y_{i}^{k}-\tilde y_{i})^\textrm{T}H_1^{\textrm{T}}}{\pi_{i}}$.}
\vspace{0.05in}
\State{{3b.} $\beta_{v_{i}}^k = \tilde \beta_{v_{i}}+\frac{(y_{i-1}-\tilde y_{i-1})^{\textrm{T}}H_2^{\textrm{T}}}{\pi_{i}}$.}
\vspace{0.05in}
\State{{3c.} $\beta_{\lambda_{i}}^k = \tilde \beta_{\lambda_{i}}-\frac{(y_{i}^k-\tilde y_{i})^{\textrm{T}}G^{\textrm{T}}}{\pi_{i}}$.}
\vspace{0.05in}
\State{{4a.} $w_i^{k} = w_{i}+\tau \left(z_{i}^{k}-\beta_{w_{i}}^k\right)$.}
\State{{4b.} $v_{i}^{k}  =  v_{i}+\tau \left(z_{i}^{k}-\beta_{v_{i}}^k\right)$.}
\State{{4c.} $\lambda_{i}^k = \lambda_{i} +\tau\left(\beta_{\lambda_{i}}^k+d-\sigma_{i}^{k}\right)$.}
\EndFor
\State{{5.} $\tilde {\mathbf w}^s = \frac{1}{T}\sum_{k=1}^T \mathbf w^k$, $\tilde {\mathbf v}^s = \frac{1}{T}\sum_{k=1}^T \mathbf v^k$,  \hfill\\\vspace{0.05in}~~~~~~~~$\tilde {\boldsymbol \lambda}^s = \frac{1}{T}\sum_{k=1}^T {\boldsymbol \lambda}^k$, and $\tilde {\boldsymbol y}^s = \frac{1}{T}\sum_{k=1}^T {\boldsymbol y}^k$ .}
 \EndWhile
 \caption{SVR-AMA for Problem~\eqref{eq:splitting_problem_condensed}.} 
 \label{alg:splittingAMA}
 \end{algorithmic} 
\end{algorithm}
\ifpaper
\begin{proof}
We summarize the logic of the proof. The inequality can be derived by the results of Theorem 5.3 in~\cite{Pu2014} by noticing that the primal updates in 
the inner loop are the same as AMA. Then, we have to take into account for Algorithm~\ref{alg:splittingAMA} that the primal variables are stochastic variables 
and that we must consider their expected values. These observations combined with the results of Theorem~\ref{th:svg_ama} lead to~\eqref{eq:th3}.
\end{proof}
\else
\begin{proof}
Refer to~\cite{ExtendedVersion}.
\end{proof}
\fi
\begin{rem}
The initial value of the dual variables $\tilde {\boldsymbol \mu}^0$ should be a feasible starting point in order to use the results of Theorem 5.3. 
This can be accomplished by noticing the following. Concerning the $\tilde \lambda_t^0$ components of $\tilde {\boldsymbol \mu}^0$, they must be in $\mathbb C_t$. 
Concerning the $\tilde w_t^0$ and $\tilde v_t^0$ components of $\tilde {\boldsymbol \mu}^0$, by providing an initial primal solution satisfying the consensus constraints 
(e.g., by using the evolution of the state starting from $x_{\textrm{init}}$ under the associated unconstrained LQR control law $u_t = K_{\textrm{LQR}}x_t$), 
they can be set equal to zero.
\end{rem}
\else
\fi
The decomposition along the length of the prediction horizon offers several advantages. First, the size of the subproblems~\eqref{eq:subproblem_t_primal} is fixed and 
independent from the length of the prediction horizon. Second, the resultant subproblems have improved numerical properties (in terms of condition number, for example), 
compared to solving Problem~\eqref{eq:initial_MPC_problem}. Third, this decomposition allows one to fully parallelize the solution of Problem~\eqref{eq:initial_MPC_problem}, 
thanks to the introduction of the consensus variables. In theory, if $N+1$ independent workers are available, the dual update of each subproblem can be assigned to its
 dedicated worker that exchanges information with its neighbors only at dedicated synchronization points to update the consensus variables, as detailed in~\cite{StathopoulosECC13}.
  If the prediction horizon, however, is larger than the number of available workers the computation of the solution has to be partially (or fully, if only one worker is available)
   serialized. This scenario can be quite common for embedded legacy control systems, where the \emph{serial} hardware architecture is formally verified and the costs to upgrade 
   to the \emph{parallel} one are too high. In this scenario, Algorithm~\ref{alg:svrama} plays a fundamental role to compute a suboptimal solution of 
   Problem~\eqref{eq:initial_MPC_problem}. 

Algorithm~\ref{alg:svrama} applied to Problem~\eqref{eq:time_splitting_MPC_problem} translates into the possibility of asynchronous updates of the independent subproblems. 
Compared to solving the subprolems in a serialized fashion (i.e., one after the other) in a synchronous framework, the asynchronous updates lead to less costly (in terms of
 computation time) iterations of the algorithm. In particular, assuming that only one worker is available, at each inner-loop iteration (steps 1-4 of the algorithm), only one 
 subproblem is randomly selected for the update. In a synchronous framework, the update of all the subproblems would have been required, which can be costly if the length of the
  horizon is large. 

Compared to other asynchronous dual algorithms (e.g.,~\cite{Notarnicola2015}), Algorithm~\ref{alg:svrama} allows one to tune and adapt (online) the probability distribution $\Pi$. 
This is particularly useful, for example, to give priority in the update to those subproblems whose associated dual variables vary the most between two iterations of the algorithm, 
as shown in the next section.   
\section{NUMERICAL EXAMPLE}
\label{sec:sec_simulation_results}
This section considers the linearized model (at a given trim condition) of an
Airbus passenger aircraft to test the proposed design. Aerospace applications offer several challenges for MPC from the computational perspective. First, these applications usually have strict real-time requirements. Second, the states have different magnitudes (ranging from few degrees to thousands of feet) affecting the conditioning of the MPC problem. Third, some of the open-loop system eigenvalues are, in general, complex conjugate close to the imaginary axis. Finally, the problem size is relatively large and a long prediction horizon is required.

We focus on the longitudinal control of the aircraft. In this respect, the
model we consider has $n=6$ states (to describe the longitudinal dynamics) and
$m=4$ control actuators.
In particular, the states associated with the longitudinal dynamics are pitch
rate [deg/sec], roll rate [deg/sec], ground speed [knots], angle of attack [deg], pitch angle [deg], and altitude [ft]. Finally, the control
surfaces for the longitudinal dynamics are the four elevators on the tail of the
aircraft.

The sampling time of the system is $T_s = 0.04$ sec and we consider an horizon length $N=60$. The total number of decision variables is 600. Furthermore, we have 3000 inequality constraints and 600 equality constraints. The goal of the MPC controller is to regulate the state of the system to the origin starting from a nonzero initial condition close to the saturation limits of the system. 

We compared the behavior of Algorithm~\ifpaper\ref{alg:splittingAMA}\fi~to a synchronous one (a state-of-the-art AMA~\cite{Tseng91,Goldstein2014}). The baseline for the comparison is the trajectory obtained using the MPC functions of MPT3~\cite{MPT3}. In particular, we are interested in showing that the possibility of tuning the probability distribution that the algorithm offers can lead to significant improvements in terms of performance of the MPC controller, especially when the solver runs for a limited number of iterations to reach a medium accuracy. In this respect, we consider four different probability distributions\ifpaper~(depicted in Figure~\ref{fig:dist})\else\fi: (i) uniform, (ii) Poisson, (iii) generalized Pareto, and (iv) adaptive. Scenario (iv) is obtained as follows. We initialize the $\Pi$ to be the uniform distribution. Then, every $T$ iterations, we check, for each $t=0,\ldots,N$ whether the following condition is verified $\|\lambda^T-\lambda^{T-1}\|^2<0.01$. If the condition is verified, $\pi_t\leftarrow 0.5\pi_t$ and the probabilities of its neighbors become $\pi_{t+1}\leftarrow\pi_{t+1}+0.25\pi_t$ and $\pi_{t-1}\leftarrow\pi_{t-1}+0.25\pi_t$.   

\ifpaper
\begin{figure}[t]
\centering 
\includegraphics[width=\linewidth]{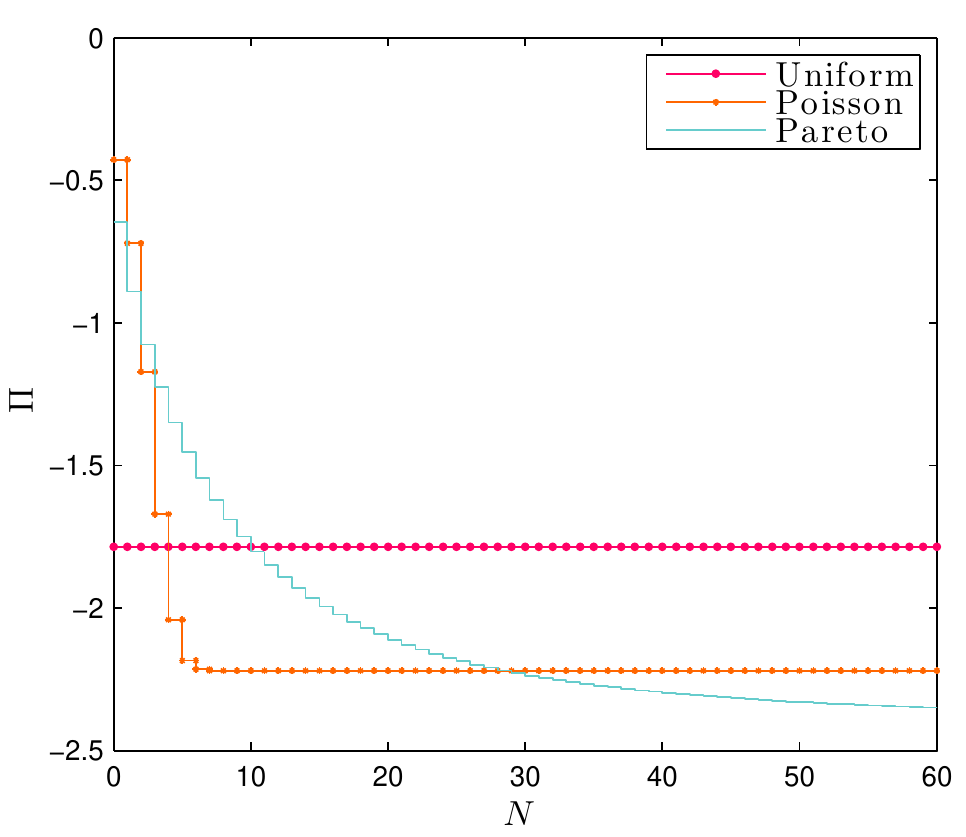}
  \caption{Probability distributions used to test the performance of Algorithm~\ref{alg:splittingAMA} plotted in $\log_{10}$ scale.}\label{fig:dist}
\end{figure}
\else
\fi
To compare the different probability tuning strategies in practice, we use the same batch size $T=10<N$, step size $\tau = 0.9$, and number of outer iterations $\bar s = 15,000$ for all the numerical experiments. From the practical point of view, we noticed that the asynchronous algorithm would allow one to select a larger step size (compared to the synchronous algorithm), despite the large condition number of the problem ($\kappa\approx 10^5$). The tuning of the parameters is made to compare with AMA. In particular, we selected the same $\tau$ and the number of iteration of AMA is set to match the number of iterations of SVR-AMA ($\bar s = 15,000\times T/ N$). This choice is made to compare, given the same computation time, the quality of the solution obtained using the synchronous and the asynchronous algorithms. 

\begin{figure*}[t] 
\minipage{0.5\textwidth}
  \includegraphics[width=1\linewidth]{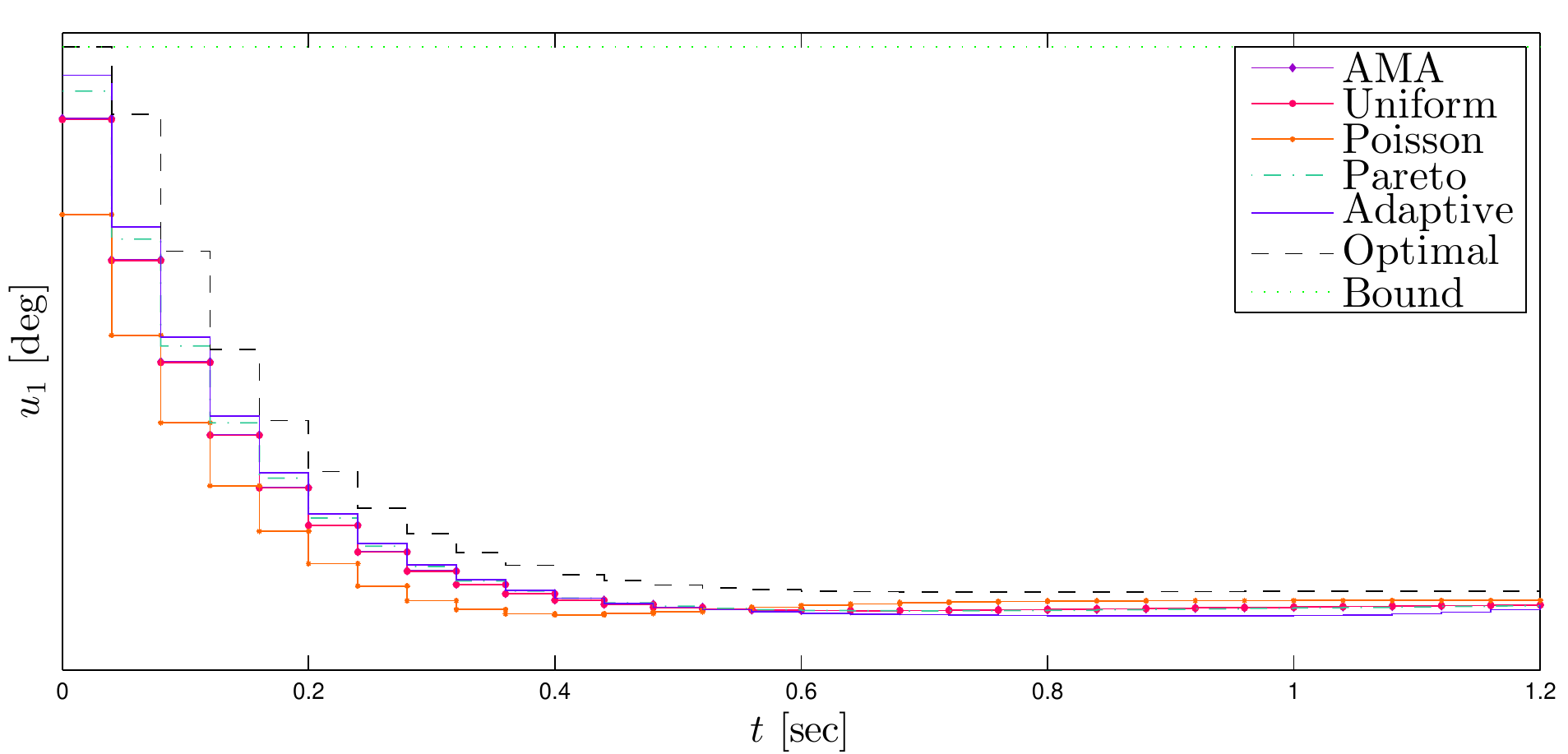}
  \caption{Control trajectories obtained using different probability distributions $\Pi$ in Algorithm~\ifpaper\ref{alg:splittingAMA}\fi~in open loop.}\label{fig:comparison_distributions_input_ol}
\endminipage~ 
\minipage{0.5\textwidth}
  \includegraphics[width=1\linewidth]{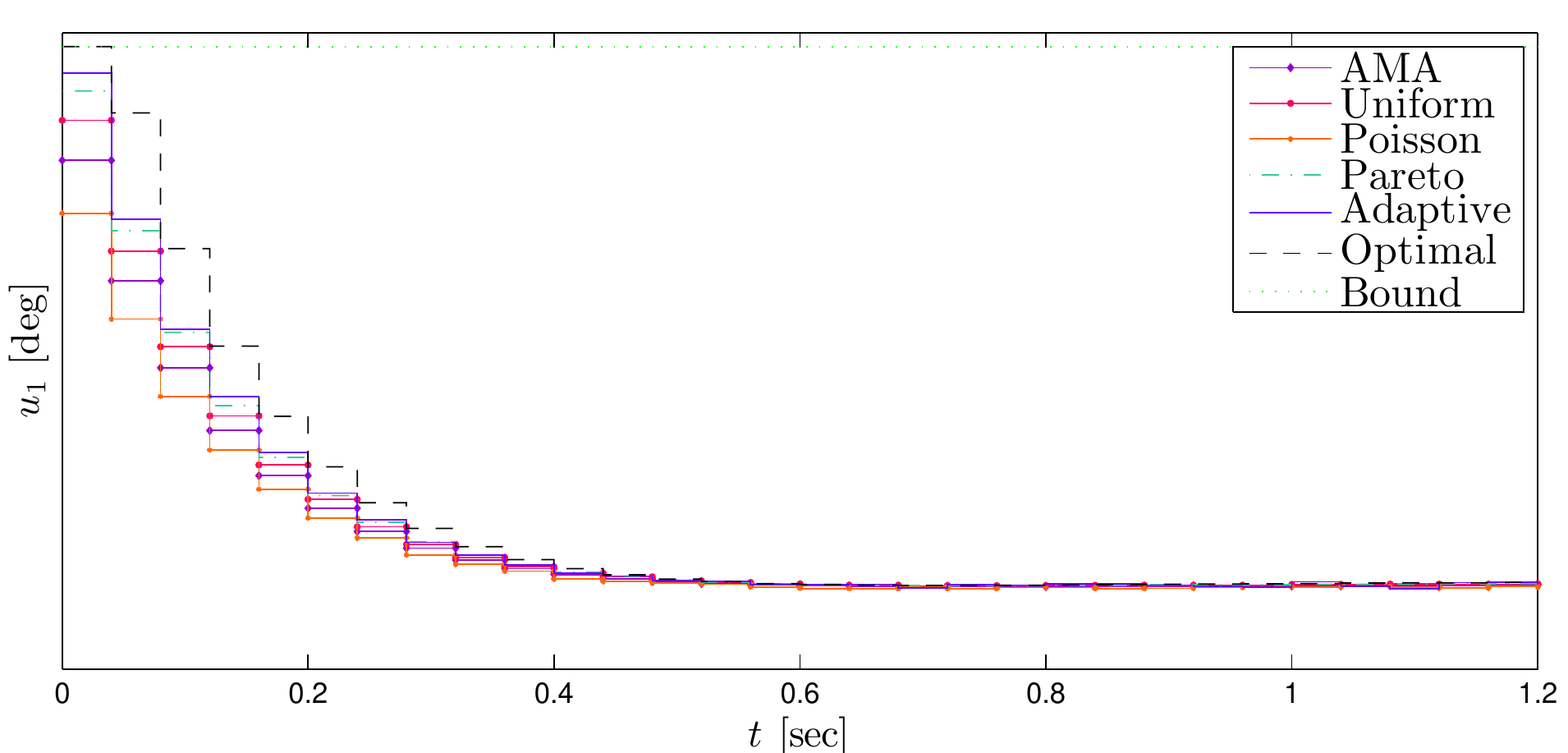}
  \caption{Control trajectories obtained using different probability distributions $\Pi$ in Algorithm~\ifpaper\ref{alg:splittingAMA}\fi~in closed loop.}\label{fig:comparison_distributions_input_cl}
\endminipage\\
\minipage{0.5\textwidth}
  \includegraphics[width=1\linewidth]{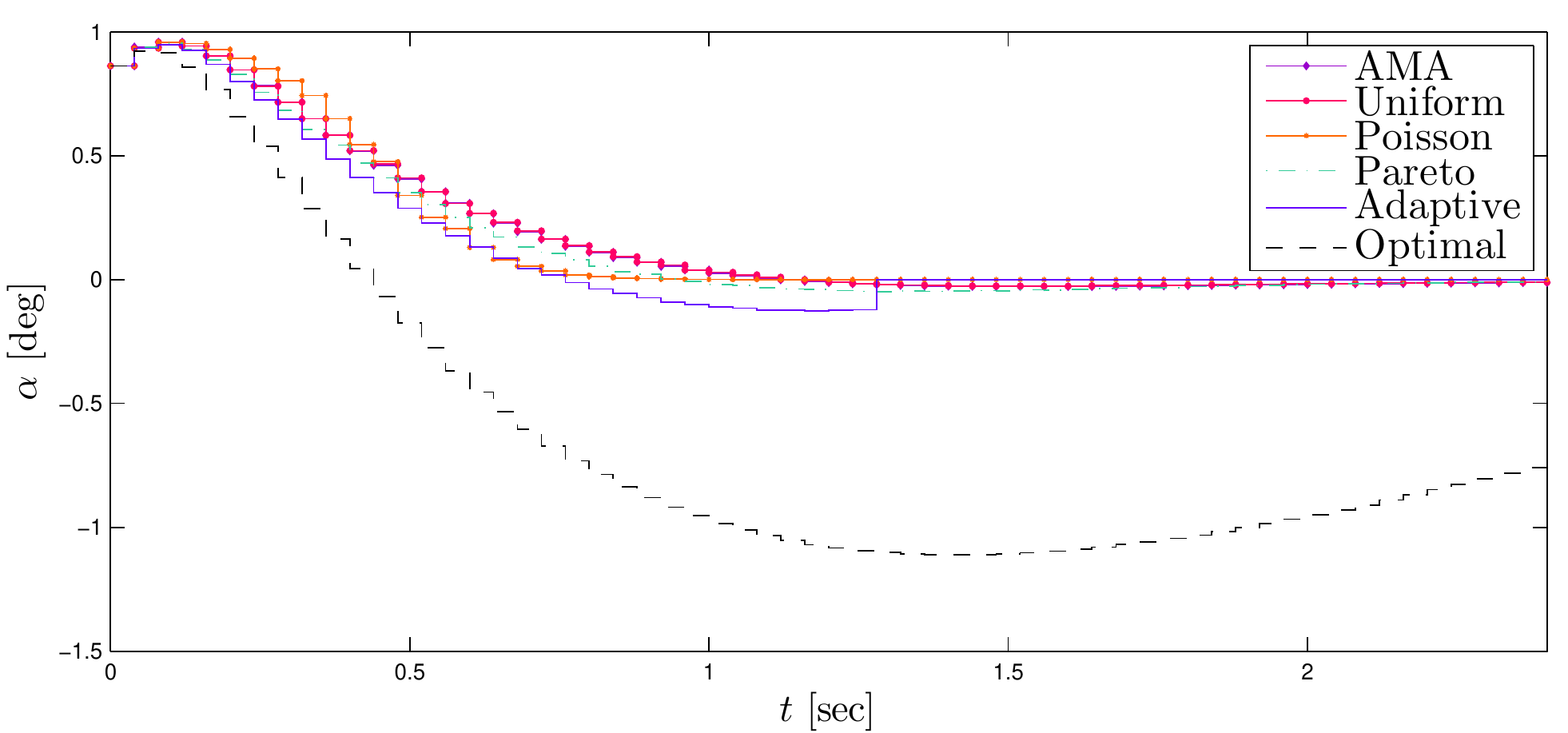}
  \caption{Angle of attack trajectories obtained using different probability distributions $\Pi$ in Algorithm~\ifpaper\ref{alg:splittingAMA}\fi~in open loop.} \label{fig:comparison_distributions_state_ol}
\endminipage~
\minipage{0.5\textwidth}
  \includegraphics[width=1\linewidth]{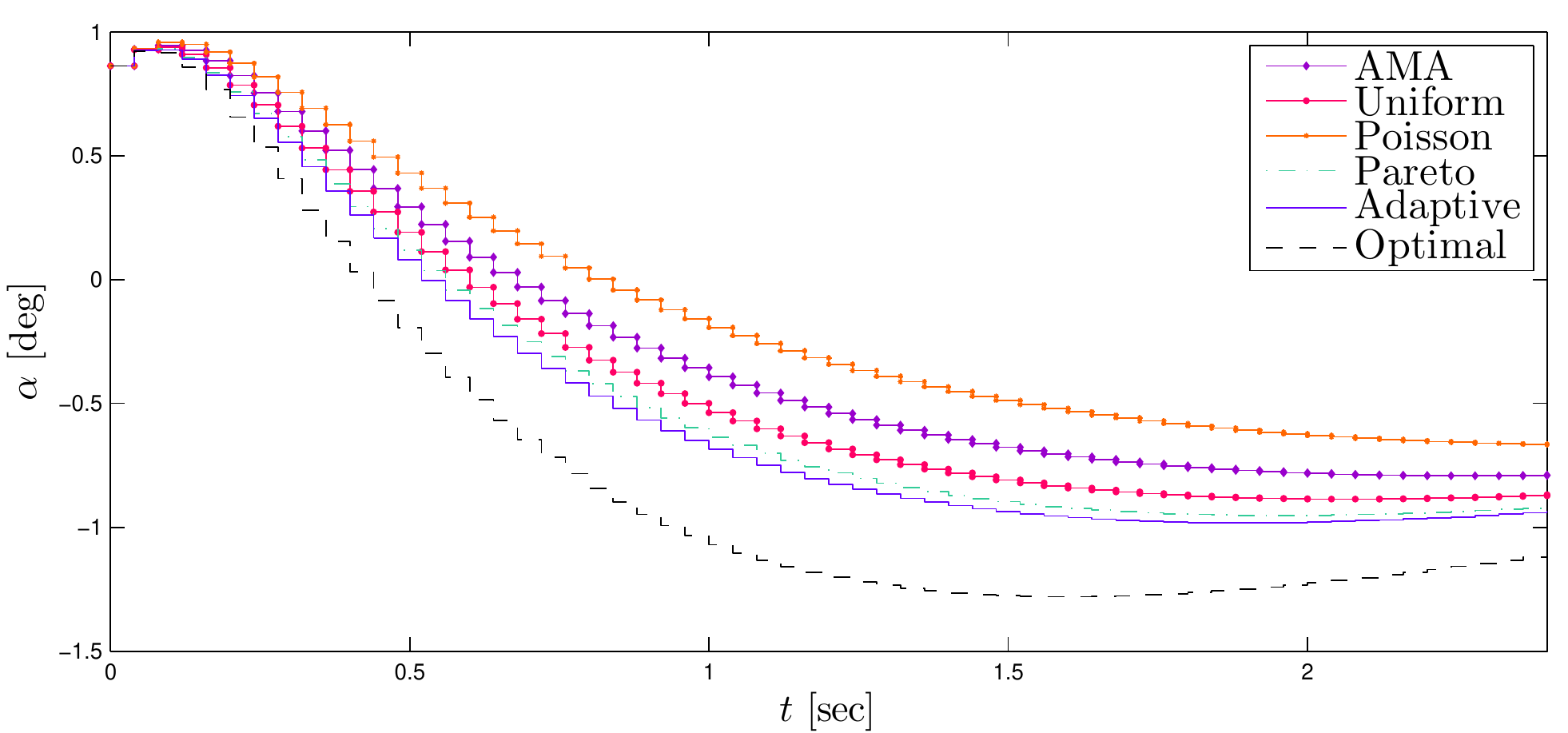}
  \caption{Angle of attack trajectories obtained using different probability distributions $\Pi$ in Algorithm~\ifpaper\ref{alg:splittingAMA}\fi~in closed loop.}\label{fig:comparison_distributions_state_cl}
\endminipage
\end{figure*}

Figures~\ref{fig:comparison_distributions_input_ol}
and~\ref{fig:comparison_distributions_state_ol} show the open-loop predictions
obtained for the elevator command and the angle of attack using different
probability distributions. As expected, the predictions are suboptimal.
Nevertheless, it is interesting to note that the results obtained using the
adaptive distribution, given the same number of iterations and step size of the
other setups. The beginning of the horizon is solved with a higher accuracy,
while the tail of the horizon is almost never updated. This does not affect the
closed-loop performance, as Figures~\ref{fig:comparison_distributions_input_cl}
and~\ref{fig:comparison_distributions_state_cl} show. In particular, the
behavior obtained using the adaptive distribution outperforms the other setups,
given the same number of iterations. In particular, compared to AMA and SVR-AMA
with uniform distribution, within the same number of iterations, SVR-AMA with adaptive distribution reaches a higher level of suboptimality, which can be a useful feature if the system has hard real-time constraints and the number of iterations of the solver allowed within the sampling time is small.

This simple, but illustrative, example clearly shows the advantages of tuning the probability distribution of the random updates, in order to improve the quality of the estimates and speed up the computation of the MPC problem solution. 

\section{CONCLUSIONS}
\label{sec:sec_conclusions}
We presented an asynchronous alternating minimization algorithm with variance reduction scheme suitable for model predictive control (MPC) applications.
 As our numerical example showed, the proposed algorithm outperforms a state-of-the art solver (the alternating minimization algorithm) in terms of number of 
 iterations needed to reach a desired level of suboptimality (i.e., when the solver terminates after a fixed number of iterations). Furthermore, the possibility of
  tuning the probability distribution of the random updates is an additional benefit for MPC applications. In particular, compared to other state-of-the-art asynchronous 
  dual solvers that only performs random updates according to a uniform distribution, the proposed algorithm allows one to prioritize the update of the variables at the
   beginning of the prediction horizon, leading to improved behavior in closed loop, as our numerical example showed. 

As part of our future work, we aim to extend these results to different splitting strategies and to formulate its accelerated version.
\ifpaper

\else

\fi 
 
\end{document}